
\input amssym.def   
\input amssym.tex   

\def \tilde{\widetilde}   

\input epsf.tex     
\epsfverbosetrue    

\def \diagram{
 \vskip \abovedisplayskip
 \epsfysize=4cm
 \centerline{\epsfbox{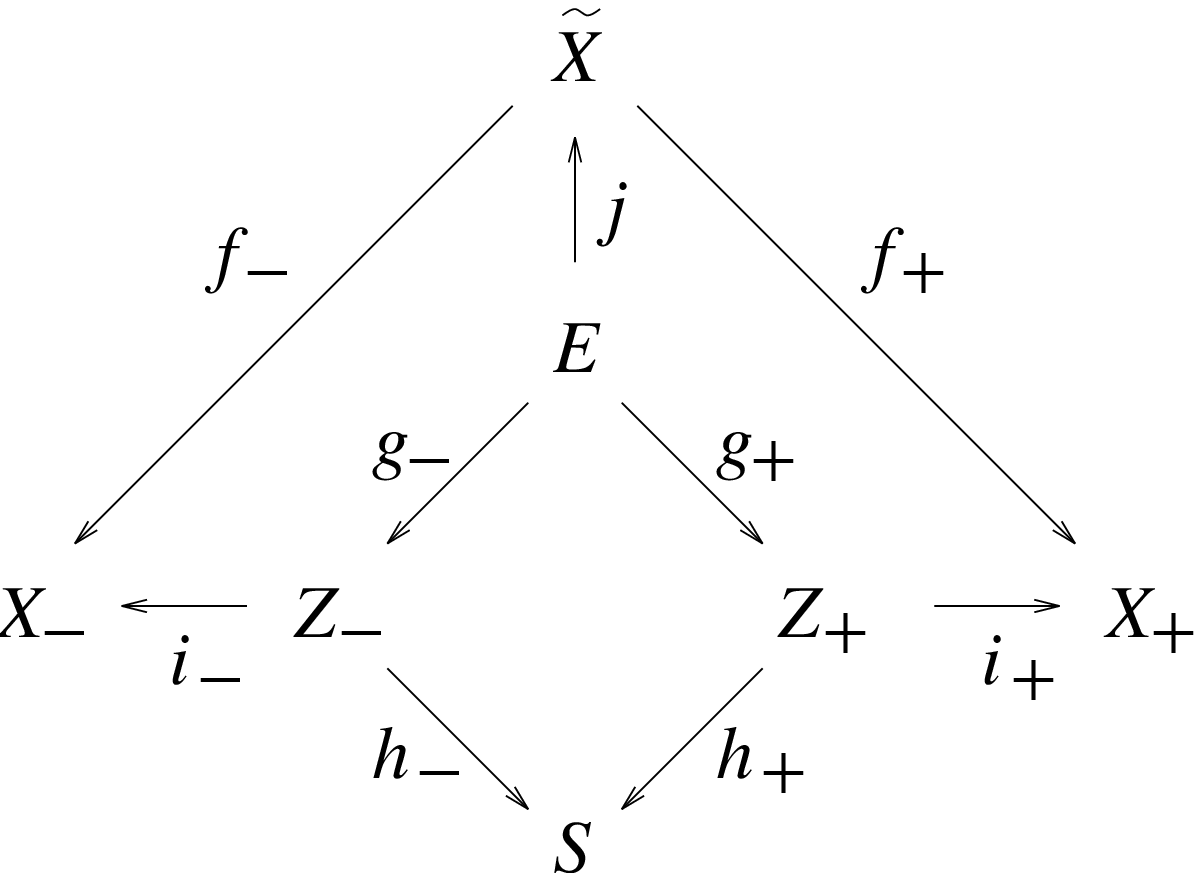}}
 \vskip -2.2cm
 \rightline{(2.1)}
 \vskip -\baselineskip
 \vskip 2.2cm
 \vskip \belowdisplayskip
 \noindent }
%
%
%


\magnification=\magstep1

\def \titlefont{\bf}
\font \authorfont=cmcsc10
\def \thanks#1{\footnote{}{{\rm #1}}}

\def \qedbox{{\vbox{\hrule\hbox{\vrule\phantom{q}\vrule}\hrule}}}
\def \qed{\qquad\qedbox\hfill\bigbreak}
\def \cite#1{{\bf [#1]}}
\def \say #1.
{\medbreak\noindent{\it#1.}\enspace}


\def \CC{{\Bbb C}}
\def \NN{{\Bbb N}}
\def \PP{{\Bbb P}}
\def \QQ{{\Bbb Q}}
\def \RR{{\Bbb R}}
\def \ZZ{{\Bbb Z}}
\def \cL{{\cal L}}
\def \cU{{\cal U}}
\def \lam{\lambda}
\def \Lam{\Lambda}
\def \sub{\subseteq}
\def \tsr{\otimes}
\def \isom{\cong}
\def \dual{^\vee}
\def \half{{1\over 2}}
\def \iff{\quad\Leftrightarrow\quad}

\def \lra{\,\longrightarrow\,}
\def \into{\hookrightarrow}
\def \birat{-\to}
\def \swa{\swarrow}
\def \sea{\searrow}
\def \bda{\big\downarrow}
\def \LRA#1{\;{\buildrel {\displaystyle #1}\over\longrightarrow}\;}
\def \BDA#1{\phantom{#1}\big\downarrow{#1}}
\def \MC{{{\cal N}_C}}
\def \JC{{J_C}}
\def \HA{H_A}
\def \PA{P_{A}}

\def \Xtil{{\tilde X}}


\ \bigskip
\centerline{\titlefont ALGEBRAIC COHOMOLOGY OF THE MODULI SPACE}
\centerline{\titlefont OF RANK 2 VECTOR BUNDLES ON A CURVE}
\bigskip
\centerline{\authorfont
V. Balaji,\thanks{Visit of the first author to Liverpool supported by
 the EU International Scientific Cooperation Initiative
 (grant CI1*-CT93-0031) and the EPSRC (grant GR/K38878).}
A.D. King\thanks{Second author supported by
 the SERC/EPSRC (grant GR/J38932).}
\& P.E. Newstead\thanks{All authors are members of
 the VBAC group of Europroj.}}
\bigskip
\centerline{16 February 1995}
\bigskip


\say Abstract.
Let $\MC$ be the moduli space of stable holomorphic vector bundles of rank 2
 and fixed determinant of odd degree, over a smooth projective curve $C$.
This paper identifies the algebraic cohomology ring $\HA^*(\MC)$,
 i.e. the subring of the rational cohomology ring $H^*(\MC;\QQ)$
 spanned by the fundamental classes of algebraic cycles,
 in terms of the algebraic cohomology ring of the Jacobian $\JC$.

\say 1991 Mathematics Subject Classification.
 \ Primary: 14H60, 14D20.\  Secondary: 14C15, 14C30, 14F25.


\beginsection 1. Introduction

Let $C$ be a smooth projective complex curve of genus $g\geq 2$
 and $\xi$ a line bundle on $C$ of odd degree $d$.
Let $\MC$ be the moduli space of rank 2 stable holomorphic vector
 bundles on $C$ with determinant $\xi$.
This notation is justified by the fact that the isomorphism class of
 the moduli space is independent of $\xi$.
In addition, let $\JC$ be the Jacobian of $C$.

For any smooth projective variety $X$,
 the (rational) algebraic cohomology group
$$
 \HA^i(X)\sub H^{2i}(X;\QQ)
$$
 is the subspace spanned by the fundamental classes of algebraic cycles
 of (complex) codimension $i$ on $X$.
(Note: unless otherwise stated, all cohomology groups in this paper,
 including Chow groups, have rational coefficients.)
The algebraic Poincar\'e polynomial of $X$ is then defined by
$$
 \PA(X;t)=\sum_{i=0}^{\infty} \dim \HA^i(X) t^i.
$$

In this paper we shall relate the algebraic cohomology of $\MC$
 to that of $\JC$.
Indeed, we shall show that $\HA^*(\MC)$ bears essentially the
 same relationship to $\HA^*(\JC)$ as the ordinary cohomology
 $H^*(\MC)$ does to $H^*(\JC)$.
This latter relationship derives from the fact that there
 is a natural isomorphism
$$
 H^*(\JC)\isom \Lam^*\bigl(H^3(\MC)\bigr).
 \eqno{(1.1)}
$$
 and that $H^*(\MC)$ is generated by
 $H^3(\MC)$ together with two other algebraic classes
 $\alpha\in H^2(\MC)$ and $\beta\in H^4(\MC)$ (see \cite{Ne2}).
Thus, there is a surjective ring homomorphism
$$
 \nu:\QQ[\alpha,\beta]\tsr H^*(\JC) \to H^*(\MC).
 \eqno{(1.2)}
$$
Furthermore (see \cite{Ha}), the ordinary Poincar\'e polynomial of $\MC$ is
$$
 P(\MC;t) = {P(\JC;t^3) - t^{2g} P(\JC;t) \over (1-t^2)(1-t^4) }.
 \eqno{(1.3)}
$$

The main point of this paper is to show that the analogues of (1.2)
 and (1.3) hold for algebraic cohomology.

\proclaim Theorem 1.
$$
 \PA(\MC;t)={\PA(\JC;t^3)-t^g\PA(\JC;t) \over (1-t)(1-t^2)}.
$$

\proclaim Theorem 2.
$$
 \HA^*(\MC) = \nu\bigl( \QQ[\alpha,\beta]\tsr \HA^*(\JC) \bigr)
$$

Note that the differences in powers of $t$ between (1.3) and Theorem 1
 are simply due to the difference in grading between ordinary
 and algebraic cohomology.
We prove Theorem 1 by using the technique of Thaddeus \cite{Th1}
 to relate a projective bundle over $\MC$ to a projective space,
 through a chain of `smooth flips' whose centres are all symmetric powers
 of the curve.
Theorem 2 then follows from Theorem 1 and the fact that
 $\nu$ takes algebraic classes (in $H^*(\JC)$) to algebraic classes.

One immediate consequence (Corollary 4.1) is that numerical and homological
 equivalence coincide for $\MC$.
In addition, the arguments of this paper may be repeated to show that
 Theorems 1 \& 2 are equally valid with $\HA^i$ replaced by
 the Hodge cohomology $H^{i,i}(X)\cap H^{2i}(X;\QQ)$.
Thus the Hodge conjecture for $\MC$ would be implied by the Hodge conjecture
 for $\JC$.
For a general curve $C$, the Hodge conjecture is known to hold for $\JC$;
 in this case, recent work of Biswas \& Narasimhan \cite{BN} shows directly
 that it also holds for $\MC$ and indeed for a large class
 of smooth moduli spaces over $C$.
On the other hand, over all curves of small genus and for moduli spaces
 of low rank bundles (i.e. $g\leq4$, $r=2$ and $g=2$, $r=3$),
 the Hodge conjecture has been verified in \cite{Bal1}.

The paper is laid out as follows.
In \S2 we describe how the algebraic cohomology
 transforms under a smooth flip.
In \S3 we prove the generalisation of Macdonald's formula (3.1)
 for the algebraic Poincar\'e polynomial of the symmetric powers of $C$.
In \S4 we prove Theorem 1 by applying the results of \S2 and \S3 to
 Thaddeus' chain of flips.
In \S5 we describe some of the consequences for a general curve.
In \S6 we deduce Theorem 2 from Theorem 1.
In \S7 we discuss how far the results can be extended to the
 Chow ring.

The work in this paper --- in particular (5.1) --- was the inspiration
 for a further investigation by the second two authors \cite{KN} into the
 structure of the ordinary cohomology ring $H^*(\MC)$.
While the two papers refer to each other to clarify various points,
 there is no strict logical dependence between them.


\say Acknowledgements.
The first author wishes to thank the Mathematics Departments at the
 University of Liverpool for their generous hospitality while this work
 was carried out,
 and also M.S. Narasimhan for very helpful discussions.
We are grateful to J. Harris for communicating to us the sketched proof
 of Proposition 5.1.


\beginsection 2. Through a flip.

In this section, we describe how the algebraic cohomology groups
 transform under the simplest type of `flip', in the sense of \cite{Th2}.
More precisely, we will say that a birational map $X_- \birat X_+$
 is a `smooth flip of type $(\lam,\mu)$ with centre $S$'
 if we have the following commutative diagram of smooth projective varieties
 \diagram
 in which the central square is Cartesian, the other two are
 blowup diagrams, $Z_+$ has codimension $\lam$ in $X_+$ and
 $Z_-$ has codimension $\mu$ in $X_-$.
Thus $g_+$ and $h_-$ are projective bundles associated
 to vector bundles of rank $\lam$ and $g_-$ and $h_+$ are projective
 bundles associated to vector bundles of rank $\mu$.
Note that the case $\lam=1$ is just that of a usual (smooth) blowup.

Recall (\cite{Fu} Prop. 6.7(e)) that for a blowup diagram
$$
 \matrix{
 E & \LRA{j} & \Xtil \cr
 \BDA{g} && \BDA{f} \cr
 Z & \LRA{i} & X }
$$
 in which $Z$ has codimension $\lam$ in $X$, the Chow groups are related
 by the fact that the map
$$
 \pmatrix{f^* & j_* \cr 0 & g_*}: A^k(X) \oplus A^{k-1}(E)
 \lra A^k(\Xtil) \oplus A^{k-\lam}(Z)
 \eqno{(2.2)}
$$
 is an isomorphism.
The analogous result is true for ordinary cohomology
 and the class map $A^k\to H^{2k}$ is natural with respect to
 both pull-back and push-forward.
Hence we also have
$$
 \HA^k(X) \oplus \HA^{k-1}(E)
 \isom \HA^k(\Xtil) \oplus \HA^{k-\lam}(Z)
$$
and thus
$$
 \PA(\Xtil) - t \PA(E) = \PA(X) - t^{\lam} \PA(Z)
$$
Applying this to both blowup diagrams in the flip diagram (2.1),
 and using the fact that $Z_+$ and $Z_-$ are projective
 bundles over $S$,  we obtain
$$
 \PA(X_+)-\PA(X_-)= {t^\lam-t^\mu \over 1-t} \PA(S)
 \eqno{(2.3)}
$$


\beginsection 3. Symmetric powers.

In this section we prove the following formula, which implicitly
 gives the algebraic Poincar\'e polynomials of $S^kC$ for all $k$.
This is a direct analogue of Macdonald's formula for the ordinary
 Poincar\'e polynomials (\cite{Mac} (4.3)).
It also turns out to be a convenient way
 to use the information.
$$
 \sum_{k=0}^{\infty} \PA(S^kC;t) s^k
 =
 {\PA(\JC;s^2t) \over (1-s)(1-st) }
 \eqno{(3.1)}
$$

{}From Collino's description of the Chow ring of $S^kC$ (\cite{Co} Theorem 3),
 one may immediately see that the algebraic cohomology ring $\HA^*(S^kC)$
 is generated by $q_k^*\bigl(\HA^*(\JC)\bigr)$,
 where $q_k:S^kC\to\JC$ is the Abel-Jacobi map,
 and the class $x\in \HA^1(S^kC)$ represented by any of the
 canonical embeddings $S^{k-1}C\into S^kC$.
Thus we have a natural surjection of rings
$$
 \phi_k:\QQ[x]\tsr\HA^*(\JC) \to \HA^*(S^kC).
$$

We can deduce (3.1) directly from the following.

\proclaim Proposition 3.1.
The restricition of $\phi_k$ to
$$
 V^k :=  \bigoplus_{{i,j\geq 0 \atop i+2j\leq k}} x^i \HA^j(\JC).
 \eqno{(3.2)}
$$
 is an isomorphism.

\say Proof.
We first prove the proposition for $k\geq 2g$,
 when $q_k:S^kC\to \JC$ is a projective bundle of
 fibre dimension $k-g$ and $x$ is its relative hyperplane class.
In this case, the restriction of $\phi_k$ to
$$
 \bigoplus_{{i,j\geq 0 \atop i\leq k-g}} x^i \HA^j(\JC).
 \eqno{(3.3)}
$$
 is an isomorphism.
The summands of (3.2) coincide with those of (3.3) when the
 degree $i+j\leq k/2$.
However, we also know that $x$ is ample (c.f. \cite{Mac}) and hence,
 by the Hard Lefschetz Theorem,
 the multiplication map $x^{k-2d}:\HA^d(S^kC)\to\HA^{k-d}(S^kC)$ is
 injective when $d\leq k/2$.
Thus $\phi_k$ is at least injective when restricted to $V^k$.
But now, $\JC$ has dimension $g$ and $\HA^*(\JC)$ satisfies
 numerical Poincar\'e duality (\cite{Lieb}).
Hence, (3.2) and (3.3) have the same dimension and
 so $\phi_k$ is also surjective when restricted to $V^k$.

For $k<2g$, the result follows by `backwards induction' based
 on the fact from \cite{Co} that
$$
 f\in\ker\phi_k \iff xf\in\ker\phi_{k+1}.
$$
More precisely, suppose that the result is true for $\phi_{k+1}$.
First observe that $xV^k\sub V^{k+1}$, so that
 $(\ker\phi_{k+1})\cap V^{k+1}=0$ implies that
 $(\ker\phi_k)\cap V^k=0$.
Secondly recall that the restriction map
 $\HA^*(S^{k+1}C)\to\HA^*(S^kC)$ is surjective,
 so that we at least know that $\phi_k(V^{k+1})=\HA^*(S^kC)$.
But now the inductive hypothesis implies that,
 if $\beta\in\HA^j(\JC)$ and $i+2j=k+1$,
 then there is a relation in $\ker\phi_{k+1}$ of the form
 $x^{i+1}\beta + \cdots$, where ``$\cdots$'' involves only
 {\it higher} powers of $x$.
Thus we may divide by $x$ to obtain a relation in $\ker\phi_k$
 of the similar form $x^i\beta + \cdots$,
 and hence $\phi_k(V^k)=\phi_k(V^{k+1})$.
\qed


\beginsection 4. Proof of Theorem 1.

We now use Thaddeus' chain of flips and the formulae from the previous
 two sections to prove Theorem 1.

Recall from \cite{Th1} the following diagram (modified slightly to
 improve the symmetry)
$$
 \matrix{
 && \Xtil_1    &&&        &&& \Xtil_w && \cr
 & \swa && \sea &&        && \swa && \sea & \cr
 X_0 &&&& X_1    & \cdots & X_{w-1} &&&& X_w \cr
 \bda        &&&&&        &&&&& \bda \cr
 pt          &&&&&        &&&&& \MC }
$$
 for the moduli space $\MC$ of bundles of degree $d=2w+1\geq 4g-3$.
Here the diagonal maps are all birational and the
 two vertical maps $X_0\to pt$ and $X_w\to\MC$ are
 projective bundles associated to vector bundles of ranks $m=d+g-1$
 and $n=d-2g+2$ respectively.
(This is, in part, a perverse way of saying that $X_0\isom \PP^{m-1}$.)
The birational map $X_{k-1}\birat X_k$ is a smooth flip of type
 $(\lam,\mu)=(k,m-2k)$ with centre $S^kC$,
 in the sense of \S2.

Thus, repeated application of (2.3) yields
$$
 \PA(X_w) = \PA(X_0) + \sum_{k=1}^w {t^k-t^{m-2k}\over 1-t}\PA(S^kC)
$$
 and then using the fact that $X_0$ and $X_w$ are projective bundles gives
$$
 (1-t^n)\PA(\MC) = \sum_{k=0}^w(t^k-t^{m-2k})\PA(S^kC)
 \eqno{(4.1)}
$$

Bringing in formula (3.1) and again exploiting the fact that $S^kC\to\JC$
 is a projective bundle for $k\geq 2g-1$, we may write, for $w\geq 2g-2$,
$$
 \eqalign{
 \sum_{k=0}^w \PA(S^kC;t)s^k &= {\PA(\JC;s^2t)\over(1-s)(1-st)}
 - \PA(\JC;t)\sum_{k=w+1}^\infty \left({1-t^{k-g+1}\over 1-t}\right) s^k \cr
 &= {\PA(\JC;s^2t)\over(1-s)(1-st)}
 - {\PA(\JC;t)\over 1-t}
 \left( {s^{w+1}\over 1-s} - {s^{w+1}t^{w-g+2}\over 1-st} \right) }
$$
Applying this to the right hand side of (4.1) yields
$$
 \eqalign{
 (1-t^n)\PA(\MC)
 &= { \PA(\JC;t^3) - \PA(\JC;t^{-3})t^{m+3} \over (1-t)(1-t^2) } \cr
 &\quad  - { \PA(\JC;t)\over 1-t }
 \left( {t^{w+1} - t^{m-g-w+1}\over 1-t}
 - {t^{2w-g+3} - t^{m-2w}\over 1-t^2} \right) }
$$
Now using Poincar\'e duality for $\JC$ to make the substitution
$$
 \PA(\JC;t^{-3})t^{3g} = \PA(\JC;t^3),
$$
 together with $m=n+3g-3$ and $2w=n+2g-3$, we obtain a factor of
 $1-t^n$ on the right hand side which cancels to leave
$$
 \PA(\MC;t)={\PA(\JC;t^3)-t^g\PA(\JC;t) \over (1-t)(1-t^2)},
 \eqno{(4.2)}
$$
 thereby proving Theorem 1.
\qed

One immediate consequence is

\proclaim Corollary 4.1.
 Numerical and homological equivalence
 coincide for $\MC$.

\say Proof.
By \cite{Lieb} Theorem 1, this is equivalent to the statement that
 Poincar\'e duality holds numerically for the algebraic cohomology
 of $\MC$, which follows from (4.2), because
 the same is true for $\JC$, by \cite{Lieb} Theorem 3.
\qed

In contrast, algebraic and homological equivalence
 do not coincide for $\MC$, for certain curves $C$ (see \cite{Bal2}).


\beginsection 5. The general curve.

We use Theorem 1 to deduce some stronger statements about $\MC$ for
the general curve, i.e. for all curves $C$ lying in the complement
of a countable union of proper closed subvarieties in the
moduli space $M_g$ of curves.

We require first the following well-known fact about the general
Jacobian.

\proclaim Proposition 5.1.
For a general curve $\HA^*(\JC)$ is generated by the $\theta$ divisor.

\say Proof.
(We give a sketch here, having not found a suitable reference.)
Any algebraic class is in $H^{p,p}$.
The subalgebra of $H^*(\JC)$ consisting of classes that are in
 $H^{p,p}$ for all curves is invariant under the monodromy action
 of the mapping class group, which factors through the obvious
 action of the symplectic group $Sp\bigl(H^1(C;\ZZ)\bigr)$.
However, the only symplectically invariant subalgebra which is
 small enough to be contained in $\bigoplus_p H^{p,p}(\JC,\CC)$ is
 the one generated by $\theta$.
On the other hand, because the Hodge filtration depends holomorphically
 on $M_g$, the condition that a given class in $H^{2p}(\JC,\QQ)$ is
 in $H^{p,p}$ determines a closed analytic subvariety of $M_g$
(or strictly Teichm\"uller space).
\qed

Hence, for a general curve $C$,
$$
 \PA(\JC) = {1-t^{g+1}\over 1-t}
$$
and thus Theorem 1 yields
$$
 \PA(\MC) = {(1-t^g)(1-t^{g+1})(1-t^{g+2})
            \over (1-t)(1-t^2)(1-t^3)}
 \eqno{(5.1)}
$$
The above proof of Proposition 5.1 actually shows that the Hodge conjecture
is true for the general Jacobian and hence, as observed in the Introduction,

\proclaim Corollary 5.2.
For general $C$, the Hodge conjecture is true for $\MC$.

In addition (5.1) can be used to deduce

\proclaim Corollary 5.3.
For a general curve, the cohomology ring $\HA^*(\MC)$ is generated
by $\alpha$, $\beta$ and $\gamma$ (in the notation of \cite{Ne2}).

\say Proof.
First note that $\alpha$, $\beta$ and $\gamma$
 are certainly algebraic classes.
{}From \cite{Ne1} and \cite{Ne2}  (c.f. \cite{KN} Proposition 2.1)
 one may see that, for $2n\leq 3g-3$, the monomials
$$
 \alpha^i\beta^j\gamma^p \qquad i+2j+3p=n, \quad i+2p < g
$$
 are independent in $\HA^n$.
Multiplying by $\alpha^{3g-3-2n}$ we may obtain an equal number of
 independent monomials in $\HA^{3g-3-n}$.
As in \cite{KN} Remark 2.3, the number of such monomials in degree $n$ is
 the coefficient of $t^n$ in
$$
 \sum_{p=0}^{[{g\over 2}]}
 {\bigl(1-t^{g-2p}\bigr)\bigl(1-t^{2g-4p}\bigr)
 \over (1-t)(1-t^2)} t^{3p}
$$
 which (c.f. \cite{KN} (2.8)) is equal to (5.1).
\qed

Alternatively, the proof of Proposition 5.1 may easily be adapted
 to prove Corollaries 5.2 and 5.3 directly.
The much harder task of extending this argument to
 all smooth moduli spaces of plain and parabolic bundles has been
 carried out in \cite{BN}.


\beginsection 6. Proof of Theorem 2.

So far, we have only found the size of the algebraic cohomology
 ring $\HA^*(\MC)$ and not identified it as a subring of the
 full cohomology ring $H^*(\MC)$.
However, Theorem 1 does indicate that there may be a natural relationship
 between the algebraic cohomology of $\JC$ and that of $\MC$,
 and this turns out to be the case.
Recall from \S1 the definition of
$$
\nu:\QQ[\alpha,\beta]\tsr H^*(\JC)\to H^*(\MC).
$$

\proclaim Proposition 6.1.
The map $\nu$ takes algebraic classes on $\JC$
 to algebraic classes on $\MC$, i.e.
$$
\nu\bigl( 1\tsr \HA^*(\JC) \bigr) \sub \HA^*(\MC).
$$

\say Proof.
We start by recalling how the isomorphism (1.1) is defined.
Let $\cL$ be a universal bundle on $C\times\JC$
 and $\phi$ the (1,1) K\"unneth component of $c_1(\cL)$.
Similarly, let $\cU$ be a universal bundle on $C\times\MC$
 and $\psi$ the (1,3) K\"unneth component of $c_2(\cU)$.
Note that, while there is an ambiguity in the choice
 of $\cL$ and $\cU$, this does not affect $\phi$ and $\psi$.
Note also that $\phi$ and $\psi$ are both algebraic classes,
 because they differ from $c_1(\cL)$ and $c_2(\cU)$ respectively
 by obviously algebraic classes.

Now $\phi$ and $\psi$ induce two correspondences
$$
 \eqalign{
 H^1(C) \to H^1(\JC)&: \omega \mapsto \int_C\omega\phi \cr
 H^1(C) \to H^3(\MC)&: \omega \mapsto \int_C\omega\psi }
 \eqno{(6.1)}
$$
 which are both isomorphisms and which can then be combined
 to give the isomorphism $H^1(\JC)\isom H^3(\MC)$.
This in turn induces the isomorphism (1.1),
 which determines the map
$$
 \nu:1\tsr H^*(\JC) \to H^*(\MC).
$$
To prove the proposition it is sufficient to show that this map
 is a correspondence induced by an algebraic class on $\JC\times\MC$.
If we define
$$
 \Delta = -\half\int_C(\phi-\psi)^2,
$$
which is clearly an algebraic class, then we claim that
$$
 \nu(1\tsr\omega) = \int_\JC \omega{\Delta^g\over g!}.
 \eqno{(6.2)}
$$

To verify the claim by direct calculation, we introduce
 a basis $e_1,\ldots, e_{2g}$ for $H^1(C;\ZZ)$ and
 let $e_1\dual,\ldots, e_{2g}\dual$ be the dual basis
 with respect to the symplectic structure given by the intersection form,
 using the convention that $\int_C e_ie_i\dual = 1$.
It is most convenient to choose a symplectic basis so that
$$
 e_i\dual = \cases{e_{i+g} & $i\leq g$ \cr -e_{i-g} & $i>g$ }
$$
We may use the isomorphisms (6.1) to define bases
 $\phi_1,\ldots,\phi_{2g}$ and $\phi_1\dual,\ldots,\phi_{2g}\dual$
 of $H^1(\JC)$, and $\psi_1,\ldots,\psi_{2g}$ and
 $\psi_1\dual,\ldots,\psi_{2g}\dual$ of $H^3(\MC)$.
With respect to these bases, we have
$$
 \phi=\sum_{i=1}^{2g} e_i\dual\phi_i
 \qquad
 \psi=\sum_{i=1}^{2g} e_i\dual\psi_i
$$
and thus
$$
 \Delta=\sum_{i=1}^g \phi_i\phi_i\dual
 + \sum_{i=1}^{2g} \phi_i\dual \psi_i
 + \sum_{i=1}^g \psi_i\psi_i\dual.
$$
Observe that the factor of $-\half$ in the definition of $\Delta$ is
 required because it is
$$
 \theta=-\half\int_C\phi^2 = \sum_{i=1}^g \phi_i\phi_i\dual
$$
 which is the ample generator of $H^1(\JC;\ZZ)$ and has the key property
$$
 \int_\JC {\theta^g\over g!} = 1.
$$
It is then fairly straightforward to verify (6.2) on monomials.
\qed

To complete the proof of Theorem 2, observe that it follows from
 results in \cite{Ne1} and \cite{Ne2} (c.f. \cite{KN} Proposition 2.1 \&
 Remark 2.2) that the map $\nu$ induces an isomorphism
$$
 \bigoplus_{(i,j,k)\in S} \alpha^i\beta^j H^k(\JC) \isom H^*(\MC)
$$
 for some subset $S\sub \NN^3$.
Moreover (c.f. \cite{KN} Remark 2.3), the identity of Poincar\'e polynomials
$$
 \sum_{(i,j,k)\in S} t^{2i+4j+3k} \dim H^k(\JC)
 =
 {P(\JC;t^3) - t^{2g} P(\JC;t) \over (1-t^2)(1-t^4) }
$$
 depends only on the fact that $\dim H^k(\JC) = \dim H^{2g-k}(\JC)$,
 i.e. that $H^*(\JC)$ satisfies Poincar\'e duality numerically.

Now Proposition 6.1 implies that
 $\nu\bigl( \QQ[\alpha,\beta]\tsr \HA^*(\JC) \bigr)\sub \HA(\MC)$
 and thus $\nu$ also embeds
$$
 \bigoplus_{(i,j,2k)\in S} \alpha^i\beta^j \HA^k(\JC)
$$
 as a subspace of $\HA^*(\MC)$.
However, we may deduce as above that
$$
 \sum_{(i,j,2k)\in S} t^{i+2j+3k} \dim \HA^k(\JC)
 =
 {\PA(\JC;t^3) - t^g\PA(\JC;t) \over (1-t)(1-t^2) }
$$
 because $\HA^*(\JC)$ also satisfies Poincar\'e duality numerically
 (by \cite{Lieb} Theorem 3).
By Theorem 1, the right hand side is $\PA(\MC)$, and so
 the subspace is equal to $\HA^*(\MC)$ and Theorem 2 is proved.

\say Remark 6.2.
Proposition 3.3(iii) of \cite{KN} identifies a natural choice for $S$
 and thereby shows that $\nu$ induces an isomorphism
$$
 \HA^*(\MC) \isom
 \bigoplus_{{i+2k<g\atop j+2k<g}} \alpha^i\beta^j \HA^k(\JC).
$$


\beginsection 7. Chow groups.

In this section, we extend the results in Section 2 to describe
 how to relate the Chow groups $A^k(X_+)$ and $A^k(X_-)$,
 when $X_+$ and $X_-$ are related by a smooth flip.
We use the notation of Section 2 and in addition
 let $\xi_\pm\in A^1(Z_\pm)$ be the relative hyperplane class for $h_\pm$.

\proclaim Proposition 7.1.
Let
$$
 \eqalign{
 B^k_+ &= \bigoplus_{s=0}^{\mu-2} (\xi_+)^s \cdot (h_+)^* A^{k-\lam-s}(S)
 = \ker(h_+)_* \sub A^{k-\lam}(Z_+) \cr
 B^k_- &= \bigoplus_{r=0}^{\lam-2} (\xi_-)^r \cdot (h_-)^* A^{k-\mu-r}(S)
 = \ker(h_-)_* \sub A^{k-\mu}(Z_-) }
$$
Then $(i_\pm)_*:B^k_\pm\to A^k(X_\pm)$ is an injection and
 there is a canonical isomorphism
$$
 A^k(X_+)/(i_+)_* B^k_+ \isom A^k(X_-)/(i_-)_* B^k_-
$$

\say Proof.
First observe that we may rewrite (2.2) as
$$
 A^k(\Xtil) = f^* A^k(X) \oplus j_* B^k,
$$
where $j_*$ is injective on
$$
 B^k =  \bigoplus_{r=0}^{\lam-2} \zeta^r \cdot g^* A^{k-r-1}(Z)
 = \ker g_* \sub A^{k-1}(E)
$$
 and $\zeta\in A^1(E)$ is the relative hyperplane class of $g$.

Further recall the `key formula' (\cite{Fu} Prop. 6.7(a)) that,
 for any $z\in A^{k-\lam}(Z)$,
$$
 f^*i_*z = j_* \bigl( \gamma\cdot g^*z \bigr)
$$
 where
$$
 \gamma = \sum_{i=0}^{\lam-1} \zeta^{\lam-1-i}\cdot g^*c_i(N).
$$
 and $N$ is the normal bundle of $Z$ in $X$.

We now use the flip diagram (2.1) to write $A^k(\Xtil)$ in two different
 ways.
{}From the right hand blowup diagram and the fact that
 $\zeta_+=g_-^*\xi_-$, we have
$$
 A^k(\Xtil) = f_+^* A^k(X_+) \oplus j_*
 \bigoplus_{r=0}^{\lam-2} g_-^* \xi_-^r \cdot g_+^* A^{k-r-1}(Z_+)
$$
However, $h_+$ is also a projective bundle and
 $\gamma_-=g_+^* \xi_+^{\mu-1} + \cdots$.
Hence, writing
$$
 B_0^k = \bigoplus_{r=0}^{\lam-2}\bigoplus_{s=0}^{\mu-2}
 g_-^* \xi_-^r\cdot g_+^*\xi_+^s\cdot g_+^*h_+^* A^{k-r-s-1}(S)
 = \ker(g_+)_*\cap\ker(g_-)_* \sub A^{k-1}(E),
$$
 we have
$$
 \eqalign{
 A^k(\Xtil)
 &= f_+^* A^k(X_+)
 \oplus j_* B_0^k
 \oplus j_* \bigoplus_{r=0}^{\lam-2}g_-^* \xi_-^r\cdot \gamma_-
 \cdot g_+^*h_+^* A^{k-\mu-r}(S) \cr
 &= f_+^* A^k(X_+)
 \oplus j_* B_0^k
 \oplus  f_-^*(i_-)_* \bigoplus_{r=0}^{\lam-2} \xi_-^r\cdot
 h_-^* A^{k-\mu-r}(S) }
$$
 using the `key formula'.
Note that $(i_-)_*$ is an injection here because $j_*$ is above.

Similarly, from the left hand blowup diagram, we obtain
$$
 A^k(\Xtil) = f_-^* A^k(X_-)
 \oplus j_* B_0^k
 \oplus f_+^*(i_+)_* \bigoplus_{s=0}^{\mu-2} \xi_+^s\cdot
 h_+^* A^{k-\lam-s}(S)
$$
Dividing $A^k(\Xtil)$ by the part common to both
expressions completes the proof.
\qed

Proposition 7.1 (with the above proof) is also valid
 with ``$A$'' replaced by ``$\HA$''
 and we may regard this as an enhanced version of (2.3).

We may use this proposition to identify the first two Chow groups
$$
 \eqalign{
 A^1(\MC) &\isom \ZZ \cr
 A^2(\MC) &\isom \cases{ A^1(C) & $g=2$ \cr A^1(C)\oplus\ZZ & $g>2$ } }
$$
The first case is well-known (c.f. \cite{Ra} Prop. 3.4), because $A^1=Pic$.
The second case is closely related to the isomorphism $\JC\isom J^2(\MC)$
 proved in \cite{MN}, where
$$
 J^2(\MC)=H^3(\MC;\RR)/H^3(\MC;\ZZ)
$$
 is the Weil-Griffiths intermediate Jacobian.
Indeed, a small modification of the argument in \cite{BV} shows that
 the Abel-Jacobi map $A_H^2(\MC)\to J^2(\MC)$ is an isomorphism,
 where $A_H^*\sub A^*$ is the ideal of cycles homologically
 equivalent to 0.

It is reasonable to hope that a closer analysis would yield
 complete information about the Chow ring of $\MC$, at least modulo
 information about the Chow ring of $\JC$.


\def \paper [#1] #2 (#3) #4; #5 <#6> #7;{
\item{{\bf[#1]}} #2,\ #4,\ {\it #5}\ {\bf #6}\ (#3)\ #7 }

\def \preprint [#1] #2; #3; #4;{
\item{{\bf[#1]}} #2,\ #3,\ #4 }

\def \book [#1] #2 (#3) #4; #5;{
\item{{\bf[#1]}} #2,\ {\it #4},\ #5 (#3) }

\beginsection References.

\preprint [Bal1] V. Balaji;
The Hodge conjecture for certain moduli varieties;
preprint 1994;

\preprint [Bal2] V. Balaji;
Cycles on some moduli spaces of vector bundles;
preprint 1994;

\preprint [BV] V. Balaji \& P.A. Vishwanath;
Deformations of Picard sheaves and moduli of pairs;
to appear in {\it Duke Math J.};

\preprint [BN] I. Biswas \& M.S. Narasimhan;
Hodge classes of moduli spaces of parabolic bundles
over the general curve;
preprint 1994;

\paper [Co] A. Collino (1975)
The rational equivalence ring of symmetric products of curves;
Illinois J. Math. <19> 567--583;

\book [Fu] W. Fulton (1984)
Intersection Theory;
Springer-Verlag;

\paper [Ha] G. Harder (1970)
Eine Bemerkung zu einer Arbeit von P.E. Newstead;
J. Reine Angew. Math. <242> 16--25;

\preprint [KN] A.D. King \& P.E. Newstead;
On the cohomology ring of the moduli space
of rank 2 vector bundles on a curve;
preprint 1995;

\paper [Lieb] D.I. Liebermann (1968)
Numerical and homological equivalence of algebraic cycles
on Hodge manifolds;
Amer. J. Math. <90> 366--374;

\paper [Mac] I.G. Macdonald (1962)
Symmetric products of an algebraic curve;
Topology <1> 319--343;

\paper [MN] D. Mumford \& P.E. Newstead (1968)
Periods of moduli spaces of vector bundles on curves;
Amer. J. Math. <90> 1201--1208;

\paper [Ne1] P.E. Newstead (1967)
Topological properties of some spaces of stable bundles;
Topology <6> 241--262;

\paper [Ne2] P.E. Newstead (1972)
Characteristic classes of stable bundles of rank 2 over an
algebraic curve;
Trans. Amer. Math. Soc. <169> 337--345;

\paper [Ra] S. Ramanan (1973)
The moduli spaces of vector bundles over an algebraic curve;
Math. Ann. <200> 69--84;

\paper [Th1] M. Thaddeus (1994)
Stable pairs, linear systems and the Verlinde formula;
Invent. Math. <117> 317--353;

\preprint [Th2] M. Thaddeus;
Geometric invariant theory and flips;
to appear in {\it J. Amer. Math. Soc.};


\begingroup
\parindent=0pt
\obeylines
\bigskip
VB:
School of Mathematics,
SPIC Science Foundation,
92 G.N. Chetty Road,
T. Nagar,
Madras 600 017, INDIA.
e-mail: balaji@ssf.ernet.in
\bigskip
ADK \& PEN:
Department of Pure Mathematics,
University of Liverpool,
P.O. Box 147,
Liverpool L69 3BX, U.K.
e-mail: aking@liv.ac.uk, newstead@liv.ac.uk
\endgroup


\bye